\documentclass[final]{article}
\usepackage{graphicx, amsmath, amssymb, amsfonts}		

\newenvironment{proof}{\begin{trivlist} \item[]{\em Proof.}}{\end{trivlist}}

\newcount\refno
\newcommand\be{\begin{equation}}
\newcommand\ee{\end{equation}}
\newcommand\bn{\begin{eqnarray}}
\newcommand\en{\end{eqnarray}}
\newcommand\bns{\begin{eqnarray*}}
\newcommand\ens{\end{eqnarray*}}

\def\eop{\hfill\rule{2.0mm}{2.0mm}}

\title{A Note on the Daubechies Approach in the Construction of Spline Type Orthogonal Scaling Functions}

\author{ Tian-Xiao He and Tung Nguyen \\
{\small Department of Mathematics}\\
{\small Student Research Honor}\\
{\small Illinois Wesleyan University}\\
{\small Bloomington, IL 61702-2900, USA}\\
}

\date{}

\begin{document}

\maketitle
\setcounter{page}{1}
\pagestyle{myheadings}
\markboth{ Spline Type Orthogonal Scaling Functions and Wavelets} 
{ T.-X. He and T. Nguyen 
}

\begin{abstract}
\noindent
{\it We use Lorentz polynomials to present the solutions explicitly of equations (6.1.7) of \cite{Dau} and (4.9) of \cite{Dau88} sot that we give an efficient way to prove Daubechies' results on the existence of spline type orthogonal scaling functions and to evaluate Daubechies scaling functions. }

\end{abstract}

\section{Introduction}
\setcounter{equation}{0}

It is well-known that Multiresolution Analysis (MRA) is a general procedure to construct wavelet basis. 

\textbf{Definition 1.1} 
\textit{
A Multiresolution Analysis (MRA) generated by function $\phi$ consists of a sequence of closed subspaces $V_j, j\in \mathbb{Z}$, of $L_2(\mathbb{R})$ satisfying
\begin{description}
\item[(i)] (nested) $V_j \subset V_{j+1}$ for all $j \in \mathbb{Z}$;
\item[(ii)] (density) $\overline{\cup_{j \in \mathbb{Z}}V_j}=L_2(\mathbb{R})$;
\item[(iii)] (separation) $\cap_{j \in \mathbb{Z}}V_j=\{0\}$;
\item[(iv)] (scaling) $f(x) \in V_j$ if and only if  $f(2x)\in V_{j+1}$ for all $j\in \mathbb{Z}$;
\item[(v)] (Basis) There exists a function $\phi \in V_0$ such that $\{\phi(x-k):k \in \mathbb{Z}\}$ is an orthonormal basis or a Riesz basis for $V_0$.\\
\end{description} 
The function whose existence asserted in (v) is called a scaling function of the MRA.
}\\\\
A scaling function $\phi$ must be a function in $L_2(\mathbb{R})$ with $\int\phi \neq 0$. Also, since $\phi \in V_0$ is also in $V_1$ and $\{\phi_{1,k}:=2^{j/2}\phi(2x-k):k\in \mathbb{Z}\}$ is a Riesz basis of $V_1$, there exists a unique sequence $\{p_k\}_{k=-\infty}^{\infty} \in l_2(\mathbb{Z})$ that describe the two-scale relation of the scaling function\\\\
\begin{equation}
\phi(x)=\sum_{k=-\infty}^{\infty} p_k \phi(2x-k),
\end{equation}
i.e., $\phi$ is of a two-scale refinable property. 
By taking a Fourier transformation on both sides of (1.1) and denoting the Fourier transformation of $\phi$ by $\hat{\phi}(\xi):=\int_{-\infty}^{\infty} \phi(x)e^{-i\xi x}dx$, we have\\\\
\begin{equation}
\hat{\phi}(\xi)=P(z)\hat{\phi}(\frac{\xi}{2}),
\end{equation}
where\\\\
\begin{equation}
P(z)=\frac{1}{2} \sum_{k=-\infty}^{\infty} p_k z^k \ and\  z=e^{-i\xi/2}
\end{equation}
Here, $P(z)$ is called the mask of the scaling function. Now, regarding the property that $\{\phi(x-k)\}$ must be an orthonormal basis, we have the following characterization theorem (see, for example, Chs. 2, 5 and 7 of Chui \cite{Chui} and Ch. 3 of Hern\'andes and Weiss \cite{HW})\\\\

\textbf{Theorem 1.2}
\textit{
Suppose the function $\phi$ satisfies the refinement relation $\phi(x)=\sum_{-\infty}^{\infty} p_k \phi(2x-k)$. Then 
(i) $\{ \phi(x-k):k\in{mathbb Z}\}$ forms an orthonormal basis only if $|P(z)|^2+|P(-z)|^2=1$ for $z\in \mathbb{C}$ with $|z|=1$.\\
(ii) Suppose $P(z)$ satisfies
\begin{enumerate} 
\item [1.] $P(z) \in C^1$ and is $2\pi$-periodic
\item [2.] $|P(z)|^2+|P(-z)|^2=1$
\item [3.] $P(1)=1$
\item [4.] $P(z) \neq 0$ for all $\xi \in [-\pi,\pi]$
\end{enumerate}
Then $\{ \phi(x-k):k\in{mathbb Z}\}$ forms an orthonormal basis.
}\\\\


Denote the cardinal B-splines with integer knots in ${\mathbb{N}}_0$ by $B_n(x)$. It is well-known that $B_n(x)$ satisfy refinement relation 

\begin{equation}
B_n(x)=\sum_{j=0}^n \frac{1}{2^{n-1}} \left(\begin{array}{c}n\\j\end{array}\right)B_n(2x-j)
\end{equation} \\\\
and have masks $P_n(z)$ such that $\hat B_n(\xi)=P_{n}(z)\hat B_{n}(\xi/2)$, where $z=e^{-i\xi /2}$ and 

\begin{equation}\label{0-1}
P_n(z)=\frac{1}{2}\sum_{j=0}^n \frac{1}{2^{n-1}}\left(\begin{array}{c}n\\j\end{array}\right)z^j=\frac{(1+z)^n}{2^n}=\left(\frac{1+z}{2}\right)^n
\end{equation}

It is clear that 
\begin{align*}
|P_n(z)|^2+|P_n(-z)|^2 &=\left|\frac{1+z}{2}\right|^{2n}+\left|\frac{1-z}{2}\right|^{2n}\\
&=cos^{2n}(\xi/4)+sin^{2n}(\xi/4) \leq cos^{2}(\xi/4)+sin^{2}(\xi/4)=1
\end{align*}
The equality happens only when n=1. Therefore,  except for the case of order one (i.e., $n=1$), $B_n(x)$ are generally not orthogonal (indeed they are Riesz basis). To induce orthogonality, Daubechies (see \cite{Dau, Dau88}) introduces a class of polynomial function factors $S(z)$. Hence, instead of $B_n(x)$, a scaling function $\phi_n(x)$, called spline type scaling functions, with the mask $P_n(z)S_n(z)$ is considered so that 

\begin{equation}
\phi_n(\xi)=P_n(z)S_n(z)\phi_n(\xi/2), 
\end{equation}
where $P_n(z)$ are defined as \eqref{0-1}. We need to construct $S_n(z)$ such that the shift set of the new scaling function form an orthogonal basis.  In other words, we need that $S_n(z)$ satisfy the following condition

\begin{equation}\label{0-2}
|P_n(z)S_n(z)|^2+|P_n(-z)S_n(-z)|^2=1
\end{equation} 
Now we consider $S_n(z)$ of the following type: $S_n(z)=a_1z+a_2z^2+...+a_nz^n$, $n\in{\mathbb N}$ and $a_i \in \mathbb{R}$, $i = 1 .. n$. When $z=1$, from equation \eqref{0-2} we have
\begin{align*}
1&=|P_n(1)S_n(1)|^2+|P_n(-1)S_n(-1)|^2\\
&=|P_n(1)|^2|S_n(1)|^2+|P_n(-1)|^2|S_n(-1)|^2\\
&=|S_n(1)|^2+0=|S_n(1)|^2
\end{align*}

Thus $S_n(1)=\sum_{i=1}^na_i=\pm 1$. From Theorem 1.2 (ii), we further impose a restriction that $\sum_{i=1}^na_i=1$ in order to ensure the orthogonality of the scaling function.\\\\
Next, we set out to find the expressions and constructions of $S_n$. We have the following Lemma and leave the proof for next section.\\\\

\textbf{Lemma 1.3}\\
\textit{ Let $S_n(z)$ be defined as above. Then there holds 
\begin{align*}
&|S_n(z)|^2 \\
&= \sum_{i=1}^n a_i^2+ 2\sum_{i=1}^{n-1}a_i a_{i+1}cos(\xi/2)+2\sum_{i=1}^{n-2}a_ia_{i+2}cos(2\xi/2)+...\\ 
&\quad +2a_1a_ncos((n-1)\xi/2).
\end{align*}
}

From Lemma 1.3, if we write each $cos(k\xi/2)$ as a polynomial of $cos(\xi/2)$, then $|S_n(z)|^2=Q_n(x)$ where $x=cos(\xi/2)$. Obviously, $Q_n(x)$ has the degree of $n-1$. It is also easy to observe that $|S_n(-z)|^2=Q_n(-x)$. Now equation (2.5) becomes
\begin{align*}
1&=|P_n(z)S_n(z)|^2+|P_n(-z)S_n(-z)|^2\\
&=cos^{2n}(\xi/4)Q_n(x)+sin^{2n}(\xi/4)Q_n(-x)\\
&=\left(\frac{1+cos(\xi/2)}{2}\right)^nQ_n(x)+\left(\frac{1-cos(\xi/2)}{2}\right)^nQ_n(-x)\\
&=\left(\frac{1+x}{2}\right)^nQ_n(x)+\left(\frac{1-x}{2}\right)^nQ_n(-x)
\end{align*}
So finally we get

\begin{equation}\label{0-3}
\left(\frac{1+x}{2}\right)^nQ_n(x)+\left(\frac{1-x}{2}\right)^nQ_n(-x)=1
\end{equation}
As a side note, \eqref{0-3} is equivalent to (6.1.7) of \cite{Dau} or (4.9) of \cite{Dau88}, but is of quite different form so that we may obtain a complete  different solution of the equation shown below in \eqref{0-4} by using Lorentz polynomials, which yields an efficient proof of sufficiency for the orthogonality and a mechanical and elementary way to construct scaling functions $\phi_{n}$.\\\\

Next, to show the existence of $Q(x)$ in the above equation, we make use of the Polynomial Extended Euclidean algorithm (see Cormen, Leiserson, Rivest, and Stein \cite{CLRS}).\\\\

\textbf{Lemma 1.4 Polynomial extended Euclidean algorithm}
\textit{
If a and b are two nonzero polynomials, then the extended Euclidean algorithm produces the unique pair of polynomials (s, t) such that as+bt=gcd(a,b), where $deg(s)<deg(b)-deg(gcd(a,b))$ and $deg(t)<deg(a)-deg(gcd(a,b))$. 
}\\\\
We notice that $gcd((\frac{1+x}{2})^n,(\frac{1-x}{2})^n)=1$, so by Lemma 1.4, there exists uniquely $Q(x)$ and $R(x)$ with degrees less than n such that $(\frac{1+x}{2})^nQ(x)+(\frac{1-x}{2})^nR(x)=1$. If we replace $x$ by $-x$ in the previous equation, we have $(\frac{1-x}{2})^nQ(-x)+(\frac{1+x}{2})^nR(-x)=1$. Due to the uniqueness of the algorithm, we conclude that $R(x)=Q(-x)$. So we have showed the existence of a unique $Q(x)=Q_n(x)$ satisfying equation \eqref{0-3}.\\\\

To construct $Q_n(x)$ explicitly, we use the Lorentz polynomials shown in Erd\'elyi and Szabados \cite{ES}, Lorentz \cite{Lor} and the following technique.
\begin{align*}
1&=\left(\frac{1+x}{2}+\frac{1-x}{2}\right)^{2n-1}\\
&=\sum_{i=0}^{2n-1}\left(\begin{array}{c}2n-1\\i\end{array}\right)\left(\frac{1+x}{2})^{2n-1-i}(\frac{1-x}{2}\right)^i\\
&=\left(\frac{1+x}{2}\right)^{n}\left[\sum_{i=0}^{n-1}\left(\begin{array}{c}2n-1\\i\end{array}\right)\left(\frac{1+x}{2}\right)^{n-1-i}\left(\frac{1-x}{2}\right)^i\right]\\
&+\left(\frac{1-x}{2}\right)^{n}\left[\sum_{i=0}^{n-1}\left(\begin{array}{c}2n-1\\i\end{array}\right)\left(\frac{1-x}{2}\right)^{n-1-i}\left(\frac{1+x}{2}\right)^i\right],
\end{align*}
where the polynomials presenting in the brackets are the Lorentz polynomials. We notice that the degrees of the two polynomials in the brackets are $n-1$, and because $Q_n(x)$ in equation 2.6 is unique, we can conclude that

\begin{equation}\label{0-4}
Q_n(x)=\sum_{i=0}^{n-1}\left(\begin{array}{c}2n-1\\i\end{array}\right)\left(\frac{1+x}{2}\right)^{n-1-i}\left(\frac{1-x}{2}\right)^i
\end{equation}
With the construction of $Q_n(x)$, we take a step further by showing the existence of $\sum a_i^2, \sum a_ia_{i+1}, ...$ in Lemma 1.3. \\\\
It is well-known that the set $\{1, cos(t), cos(2t), ..., cos((n-1)t)\}$ is linearly independent. As a result, $\{1, cos(\xi/2), cos(2\xi/2), ..., cos((n-1)\xi/2)\}$ forms a basis of the space $P_{n-1}(x)=\{P(x): x=cos(\xi/2)$ and P is a polynomial of degree less than $n$. 
Based on this fact and the existence of $Q_n(x)$ in equation (2.6), it is obvious that the coefficients $\sum a_i^2, \sum a_ia_{i+1}, ...$ in Lemma 2.2 must exist uniquely.

We now survey the results of \$6.1 of \cite{Dau} and \$4.B. of \cite{Dau88} in the following theorem and leave its proof by using our mask expression 
\eqref{0-4} for next section.

\textbf{Theorem 1.5}
\textit{
Let $P_{n}(z)$ and $S_{n}(z)$ be defined as above. Then for $n=1,2,\ldots$, the spline type function $\phi_{n}(x)$ with the mask $P(z)=P_{n}(z)S_{n}(z)$, where $S_n(z)=a_1z+a_2z^2+...+a_nz^n$, is a scaling function that generates an orthogonal basis of $V_{0}$ in its MRA. } 

Although we can not give the explicit expression of $\phi_{n}(x)$, we may use the recursive method presented in Theorem 5.23 of \cite{BN} by Boggess and Narcowich to find an approximation of $\phi_{n}(x)$ with any accuracy. It is easy to see that all three conditions required by Theorem 5.23 are satisfied by $P_{n}(z)$: (i) $P_{n}(1)=1$, (ii) $|P_{n}(z)|^{2}+|P_{n}(-z)|^{2}=1$ ($|z|=1$), and (iii) $P_{n}(z)|>0$ ($\xi in [-\pi, \pi]$). The examples for the cases $n=3$ and $4$ will be provided in Section $3$ to demonstrate this procedure, while all proofs are given in next section. 

\section{Proofs}

\textit{Proof of Lemma 1.3} We have
\begin{align*}
&|S_n(z)|^2\\
&=|a_1(cos(\xi/2)-isin(\xi/2))+...+a_n(cos(n\xi/2)-isin(n\xi/2))|^2\\
&=|(a_1cos(\xi/2)+...+a_ncos(n\xi/2))-i(a_1sin(\xi/2)+...+a_nsin(n\xi/2))|^2\\
&=(a_1cos(\xi/2)+...+a_ncos(n\xi/2))^2+(a_1sin(\xi/2)+...+a_nsin(n\xi/2))^2\\
&=a_1(cos^2(\xi/2)+sin^2(\xi/2))+...+a_n(cos^2(n\xi/2)+sin^2(n\xi/2))\\
&\quad +\sum_{i\neq j}2a_ia_j(cos(i\xi/2)cos(j\xi/2)+sin(i\xi/2)sin(j\xi/2))\\
\end{align*}
\begin{align*}
&= \sum_{i=1}^n a_i^2+\sum_{i\neq j}2a_ia_jcos((i-j)\xi/2)\\
&=\sum_{i=1}^n a_i^2+ 2\sum_{i=1}^{n-1}a_i a_{i+1}cos(\xi/2)+2\sum_{i=1}^{n-2}a_ia_{i+2}cos(2\xi/2)+...\\
&\quad +2a_1a_ncos((n-1)\xi/2)
\end{align*}

A similar procedure can be applied to find $|S_n(-z)|^2$
\begin{align*}
&|S_n(-z)|^2 \\
&= \sum_{i=1}^n a_i^2- 2\sum_{i=1}^{n-1}a_i a_{i+1}cos(\xi/2)+2\sum_{i=1}^{n-2}a_ia_{i+2}cos(2\xi/2)+...\\
&\quad +(-1)^n2a_1a_ncos((n-1)\xi/2).
\end{align*}

We now prove Theorem 1.5 using our new mask expression \eqref{0-4}. 
\begin{proof}
From theorem 1.2 (ii), the sufficient conditions for the orthogonality of the scaling function are
\begin{enumerate} 
\item [1.] $P(z) \in C^1$ and is $2\pi$-periodic
\item [2.] $|P(z)|^2+|P(-z)|^2=1$
\item [3.] $P(1)=1$
\item [4.] $P(z) \neq 0$ for all $\xi \in [-\pi,\pi]$
\end{enumerate}
From the construction of our $P(z)=P_n(z)S_n(z)$, the first two conditions are automatically satisfied. The third condition is also obvious: $P(1)=P_n(1)S_n(1)=\left(\frac{1+1}{2}\right)^n \sum_{i=1}^na_i=1$ according to the construction of $S_n(z)$. Now we will prove that the final condition is fulfilled as well.\\\\
Indeed, if $\xi \in [-\pi,\pi]$, then firstly we have
\begin{align*}
|P_n(z)|&=|P_n(e^{-i\xi/2})|\\
&=\left|\frac{1+e^{-i\xi/2}}{2}\right|^n=\left|\frac{1+cos(\xi/2)-isin(\xi/2)}{2}\right|^n\\
&=|cos^2(\xi/4)-isin(\xi/4)cos(\xi/4)|^n=\sqrt{cos^4(\xi/4)+cos^2(\xi/4)sin^2(\xi/4)}^n\\
&=|cos(\xi/4)|^n \geq |cos(\pi/4)|^n>0 \ for\ \xi \in [-\pi,\pi]  
\end{align*}
Secondly, from equation \eqref{0-4} we have
\begin{align*}
|S_n(z)|^2= Q_n(x)&=\sum_{i=0}^{n-1}\left(\begin{array}{c}2n-1\\i\end{array}\right)\left(\frac{1+x}{2}\right)^{n-1-i}\left(\frac{1-x}{2}\right)^i\\
&\geq \sum_{i=0}^{n-1}\left(\begin{array}{c} n-1\\i\end{array}\right)\left(\frac{1+x}{2}\right)^{n-1-i}\left(\frac{1-x}{2}\right)^i\\
&=\left(\frac{1+x}{2}+\frac{1-x}{2}\right)^{n-1}=1
\end{align*}
Thus $|S_n(z)|\geq 1$ and $|P(z)|=|P_n(z)||S_n(z)|\geq|cos(\pi/4)|^n >0$. 

The remaining thing is to show that the newly constructed scaling function $\phi_{n}(x)$ with mask $P_{n}(z)S_{n}(z)$ is in $L_2(\mathbb{Z})$, where $S_n(z)=\sum_{j=1}^n a_jz^j$ 
and $n\in {\mathbb N}$. From He \cite{He01, He04}, we know that $\phi_{n} \in L_2(\mathbb{R})$ if 
\begin{equation}
 n \sum_{j=1}^n a_j^2 < 2^{2n-1} 
\end{equation}
Recall from Lemma 1.3 that
\begin{align*}
&Q_{n}(x)=Q_{n}(cos(\xi/2))\\
&= \sum_{i=1}^n a_i^2+ 2\sum_{i=1}^{n-1}a_i a_{i+1}cos(\xi/2)+...+2a_1a_ncos((n-1)\xi/2)
\end{align*}
Taking the integration from 0 to $2\pi$ of both sides, we have\\
\begin{align*}
&\int_0^{2\pi}Q_{n}(x)d\xi\\
&=\int_0^{2\pi}\left( \sum_{i=1}^n a_i^2+ 2\sum_{i=1}^{n-1}a_i a_{i+1}cos(\xi/2)+...+2a_1a_ncos((n-1)\xi/2)\right)d\xi\\
&=2\pi\sum_{i=1}^n a_i^2
\end{align*}
On the other hands, from the expression of $Q_{n}(x)$ in \eqref{0-4} we have\\
\begin{align*}
&\int_0^{2\pi}Q_{n}(x)d\xi=\int_0^{2\pi}\sum_{i=0}^{n-1}\left(\begin{array}{c}2n-1\\i\end{array}\right)\left(\frac{1+x}{2}\right)^{n-1-i}\left(\frac{1-x}{2}\right)^i d\xi\\
&=\int_0^{2\pi}\sum_{i=0}^{n-1}\left(\begin{array}{c}2n-1\\i\end{array}\right)\left(\frac{1+cos(\xi/2)}{2}\right)^{n-1-i}\left(\frac{1-cos(\xi/2)}{2}\right)^i d\xi\\
\end{align*}
Combining the two equations above we have\\
\begin{equation}\label{1-1}
\sum_{i=1}^n a_i^2=\frac{1}{2\pi}\int_0^{2\pi}\sum_{i=0}^{n-1}\left(\begin{array}{c}2n-1\\i\end{array}\right)\left(\frac{1+cos(\xi/2)}{2}\right)^{n-1-i}\left(\frac{1-cos(\xi/2)}{2}\right)^i d\xi
\end{equation}
Now we need to show that the expression on the right hand side of \eqref{1-1} is smaller than $\frac{2^{2n-1}}{n}$. \\\\
First of all, it is easy to see that for $0\leq i\leq n-1$
\begin{align*}
\left(\begin{array}{c}2n-1\\i\end{array}\right)<\left(\begin{array}{c}2n-1\\n-1\end{array}\right)=\frac{1}{2}\left(\begin{array}{c}2n\\n\end{array}\right)
\end{align*}
Applying this inequality into \eqref{1-1} yields
\begin{align}\label{1-2}
&\frac{1}{2\pi}\sum_{i=0}^{n-1}\left(\begin{array}{c}2n-1\\i\end{array}\right)\int_0^{2\pi}\left(\frac{1+cos(\xi/2)}{2}\right)^{n-1-i}\left(\frac{1-cos(\xi/2)}{2}\right)^i d\xi\nonumber\\
&<\frac{1}{4\pi}\left(\begin{array}{c}2n\\n\end{array}\right)\sum_{i=0}^{n-1}\int_0^{2\pi}\left(cos\frac{\xi}{4}\right)^{2(n-1-i)}\left(sin\frac{\xi}{4}\right)^{2i} d\xi\nonumber\\
&=\frac{1}{\pi}\left(\begin{array}{c}2n\\n\end{array}\right)\sum_{i=0}^{n-1}\int_0^{\pi/2}(cosx)^{2(n-1-i)}(sinx)^{2i} dx
\end{align}
for $x=\xi/4$.\\\\
Now let 
\begin{align*}
A=\sum_{i=0}^{n-1}\int_0^{\pi/2}(cosx)^{2(n-1-i)}(sinx)^{2i} dx
\end{align*}
We can express A as the sumation of two terms $A=A_1+A_2$, where
\begin{align*}
&A_1=\sum_{i=0}^{[\frac{n-1}{2}]}\int_0^{\pi/2}(cosx)^{2(n-1-i)}(sinx)^{2i} dx\\
&A_2=\sum_{i=[\frac{n-1}{2}]+1}^{n-1}\int_0^{\pi/2}(cosx)^{2(n-1-i)}(sinx)^{2i} dx
\end{align*}
For $0\leq i \leq [\frac{n-1}{2}], 2(n-1-i)>2i$, the term $A_1$ becomes
\begin{align*}
A_1&=\sum_{i=0}^{[\frac{n-1}{2}]}\int_0^{\pi/2}(cosx)^{2(n-1-2i)}(sinx cosx)^{2i} dx\\
&=\sum_{i=0}^{[\frac{n-1}{2}]}\frac{1}{4^i}\int_0^{\pi/2}(cosx)^{2(n-1-2i)}(sin2x)^{2i} dx\\
&\leq  \sum_{i=0}^{[\frac{n-1}{2}]}\frac{1}{4^i}\int_0^{\pi/2}(cosx)^{2(n-1-2i)}dx
\end{align*}
For $[\frac{n-1}{2}]+1 \leq i \leq n-1, 2(n-1-i)<2i$, the term $A_2$ becomes
\begin{align*}
A_2&=\sum_{i=[\frac{n-1}{2}]+1}^{n-1}\int_0^{\pi/2}(sinxcosx)^{2(n-1-i)} (sinx)^{2(2i-n+1)} dx\\
&=\sum_{i=[\frac{n-1}{2}]+1}^{n-1}\frac{1}{4^{n-1-i}}\int_0^{\pi/2} (sinx)^{2(2i-n+1)}(sin2x)^{2(n-1-i)} dx\\
&\leq  \sum_{i=[\frac{n-1}{2}]+1}^{n-1}\frac{1}{4^{n-1-i}}\int_0^{\pi/2}(sinx)^{2(2i-n+1)}dx\\
&\leq \sum_{i=0}^{[\frac{n-1}{2}]}\frac{1}{4^i}\int_0^{\pi/2}(sinx)^{2(n-1-2i)}dx
\end{align*}
Next, we make use of the following well-known result
\begin{align*}
\int_0^{\pi/2}(sinx)^{2n}dx &=\int_0^{\pi/2}(cosx)^{2n}dx=\frac{(2n-1)!!}{(2n)!!}\frac{\pi}{2}\\
&=\frac{(2n)!}{[(2n)!!]^2}\frac{\pi}{2}=\frac{\pi}{2}\frac{(2n)!}{4^n(n!)^2}=\frac{\pi}{2}\frac{1}{4^n}\left(\begin{array}{c}2n\\n\end{array}\right)
\end{align*}
Next, we try to find the upper bound for $\left(\begin{array}{c}2n\\n\end{array}\right)$. Based on Stirling estimation in [16], we have the following inequalities
\begin{equation}\label{1-3}
\left(\begin{array}{c}2n\\n\end{array}\right) \leq \frac{4^n}{\sqrt{3n+1}}
\end{equation}
and 
\begin{equation}\label{1-4}
\left(\begin{array}{c}2n\\n\end{array}\right) \leq \frac{4^n}{\sqrt{\pi n}}\left(1+\frac{1}{12n-1}\right)
\end{equation}
Using (3.6) on $A_1$ and $A_2$ yields
\begin{align*}
A&=A_1+A_2\leq 2\sum_{i=0}^{[\frac{n-1}{2}]}\frac{1}{4^i}\frac{1}{\sqrt{3(n-1-2i)+1}} \frac{\pi}{2}\\
&=\pi \sum_{i=0}^{[\frac{n-1}{2}]}\frac{1}{3^i}\frac{1}{(\frac{4}{3})^i\sqrt{3(n-1-2i)+1}}
\end{align*}
We consider the denominator of the fraction, and let 
\begin{align*}
f(x)=(\frac{4}{3})^{2x}[3(n-1-2x)+1], \ x\in [0, \frac{n-1}{2}]
\end{align*}
Surveying the function, we have f(x) attains minimum at 0, or
\begin{align*}
f(x)=(\frac{4}{3})^{2x}[3(n-1-2x)+1]>f(0)=3(n-1)+1
\end{align*} 
Thus, we have
\begin{equation}\label{1-5}
A\leq \pi \sum_{i=0}^{[\frac{n-1}{2}]}\frac{1}{3^i}\frac{1}{\sqrt{3(n-1)+1}}\leq \frac{\pi}{\sqrt{3n-2}}\sum_{i=0}^{\infty}\frac{1}{3^i}=\frac{3\pi}{2\sqrt{3n-2}}
\end{equation}
Finally, combining \eqref{1-1}-\eqref{1-5}, 
we have
\begin{align*}
\sum_{i=1}^n a_i^2 \leq \frac{1}{\pi}\left(\begin{array}{c}2n\\n\end{array}\right) A \leq \frac{1}{\pi}\frac{4^n}{\sqrt{\pi n}}\left(1+\frac{1}{12n-1}\right)\frac{3\pi}{2\sqrt{3n-2}}
\end{align*}
For $n\geq 17$, we can easily verify that the right hand side is less than $\frac{2^{2n-1}}{n}$. Using Mathematica for direct calculation of the case $n \leq 16$, we find that the inequality in theorem 3.1 holds. Thus, it holds every interger n, and we have shown that $\phi \in L_2(\mathbb{R})$ and complete the proof.
\end{proof}
\eop\\

\section{Examples}

It is easy to find that $S_{1}(z)=z$ and $\phi_{n}(x)$ is the Haar function. For $n=2$, 

\[
S_{2}(z)=\frac{1+\sqrt{3}}{2}z+\frac{1-\sqrt{3}}{2}z^{2}
\]
and the corresponding $\phi_{2}(x)$ is the Daubechies $D_{2}$ scaling function.

As examples, we consider the efficiency of the computation in using mask expression \eqref{0-4} to construct Daubechies scaling functions (the spline type scaling functions) for the cases of $n=3$ and $4$. According to the previous sections, we will construct the Daubechies $D_{3}$ scaling function $\phi_3(x)$ from the third order B-spline function $B_3(x)$ using our expression \eqref{0-4}. \\\\
In order to construct the function $\phi_3(x)$, we start with its mask $P_3(z)S_3(z)$, where $P_3(z)=(\frac{1+z}{2})^3$ is the mask of the third order B-spline. Let $S_3(z)=a_1z+a_2z^2+a_3z^3$, then by Lemma 1.3, we have
\begin{align*}
Q_3(x)=|S_3(z)|^2&=(a_1^2+a_2^2+a_3^2)+2(a_1a_2+a_2a_3)cos(\xi/2)+2a_1a_3cos(\xi)\\
&=(a_1^2+a_2^2+a_3^2)+2(a_1a_2+a_2a_3)cos(\xi/2)+2a_1a_3(2cos^2(\xi/2)-1)\\
&=(a_1^2+a_2^2+a_3^2-2a_1a_3)+2(a_1a_2+a_2a_3)cos(\xi/2)+4a_1a_3cos^2(\xi/2)\\
&=(a_1^2+a_2^2+a_3^2-2a_1a_3)+2(a_1a_2+a_2a_3)x+4a_1a_3x^2
\end{align*}
where $x=cos(\xi/2)$ and $z=e^{-i\xi/2}$.\\\\
On the other hand, by equation \eqref{0-4} we have 
\begin{align*}
Q_3(x)&=\sum_{i=0}^{2}\left(\begin{array}{c}5\\i\end{array}\right)\left(\frac{1+x}{2}\right)^{2-i}\left(\frac{1-x}{2}\right)^i\\
&=\left(\frac{1+x}{2}\right)^2+5\left(\frac{1+x}{2}\right)\left(\frac{1-x}{2}\right)+10\left(\frac{1-x}{2}\right)^2\\
&=\frac{3}{2}x^2-\frac{9}{2}x+4
\end{align*}
Thus, from the above equations, we have the following system of equations
\begin{equation}
\begin{cases}
a_1^2+a_2^2+a_3^2-2a_1a_3  =  4\\
2(a_1a_2+a_2a_3)  =  -\frac{9}{2}\\
4a_1a_3  =  \frac{3}{2}
\end{cases}
\end{equation}
Simplify (4.1) we get
\begin{equation}
\begin{cases}
a_1^2+a_2^2+a_3^2  =  \frac{19}{4}\\
a_1a_2+a_2a_3  =  -\frac{9}{4}\\
a_1a_3  =  \frac{3}{8}
\end{cases}
\end{equation}
From this system, we have $(a_1+a_2+a_3)^2=a_1^2+a_2^2+a_3^2 +2(a_1a_2+a_2a_3+a_1a_3)=1$. Without loss of generality, consider the case $a_1+a_2+a_3=1$. Combining this with $a_1a_2+a_2a_3  =  -\frac{9}{4}$ and $a_1a_3  =  \frac{3}{8}$ we have the following solution
\begin{equation}
\begin{cases}
a_1=\frac{1}{4}\left(1+\sqrt{10}-\sqrt{5+2\sqrt{10}}\right)\\
a_2=\frac{1}{2}\left(1-\sqrt{10}\right)\\
a_3=\frac{1}{4}\left(1+\sqrt{10}+\sqrt{5+2\sqrt{10}}\right)
\end{cases}
\end{equation}
We verify the condition for $\phi_3$ to be in $L_2(\mathbb{R})$
\begin{align*}
a_1^2+a_2^2+a_3^2=\frac{19}{4} < \frac{32}{3}=\frac{2^{2 \cdot 3 -1}}{3}
\end{align*}
Thus $\phi_3$ is indeed in $L_2(\mathbb{R})$. Now we will attempt to construct $\phi_3$ explicitly. It has the mask
\begin{align*}
P_3(z)S_3(z)&=\left(\frac{1+z}{2}\right)^3 (a_1z+a_2z^2+a_3z^3)\\
&=0.0249 x - 0.0604 x^2 - 0.095 x^3 + 0.325 x^4 + 0.571 x^5 + 0.2352 x^6
\end{align*}

Hence, 
we have the refinement equation
\begin{equation}
\begin{split}
\phi_3(x)&=0.0498 \phi_3(2x-1) - 0.121 \phi_3(2x-2) - 0.191 \phi_3(2x-3) \\
&+ 0.650 \phi_3(2x-4) + 1.141 \phi_3(2x-5) + 0.4705 \phi_3(2x-6)
\end{split}
\end{equation}

Again, we construct the Daubechies $D_{4}$ orthogonal scaling function $\phi_4(x)$ from the forth order B-spine $B_4(x)$ with the mask $P_4(z)= (\frac{1+z}{2})^4$. \\\\
Now we examine the mask $P_4(z)S_4(z)$ of $\phi_4(x)$ where we define $S_4(z)=a_1z+a_2z^2+a_3z^3+a_4z^4$. By Lemma 1.3, we have
\begin{align*}
Q_4(x)=|S_4(z)|^2&=(a_1^2+a_2^2+a_3^2+a_4^2)+2(a_1a_2+a_2a_3+a_3a_4)cos(\xi/2)\\
&+2(a_1a_3+a_2a_4)cos(\xi)+2a_1a_4cos(3\xi/2)\\
&=(a_1^2+a_2^2+a_3^2+a_4^2)+2(a_1a_2+a_2a_3+a_3a_4)cos(\xi/2)\\
&+2(a_1a_3+a_2a_4)(2cos^2(\xi/2)-1)+2a_1a_4(4cos^3(\xi/2)-3cos(\xi/2))\\
&=(a_1^2+a_2^2+a_3^2+a_4^2-2a_1a_3-2a_2a_4)\\
&+(2a_1a_2+2a_2a_3+2a_3a_4-6a_1a_4)cos(\xi/2)\\
&+(4a_1a_3+4a_2a_4)cos^2(\xi/2)+8a_1a_4cos^3(\xi/2)\\
&=(a_1^2+a_2^2+a_3^2+a_4^2-2a_1a_3-2a_2a_4)\\
&+(2a_1a_2+2a_2a_3+2a_3a_4-6a_1a_4)x\\
&+(4a_1a_3+4a_2a_4)x^2+8a_1a_4x^3
\end{align*}
where $x=cos(\xi/2)$ and $z=e^{-i\xi/2}$.\\\\
Using equation \eqref{0-4} we get another expression for $Q_4(x)$
\begin{align*}
Q_4(x)&=\sum_{i=0}^{3}\left(\begin{array}{c}7\\i\end{array}\right)\left(\frac{1+x}{2}\right)^{3-i}\left(\frac{1-x}{2}\right)^i\\
&=\left(\frac{1+x}{2}\right)^3+7\left(\frac{1+x}{2}\right)^2\left(\frac{1-x}{2}\right)+21\left(\frac{1+x}{2}\right)\left(\frac{1-x}{2}\right)^2+35\left(\frac{1-x}{2}\right)^3\\
&=8-\frac{29}{2}x+10x^2-\frac{5}{2}x^3
\end{align*}
Thus, from the above equations, we have the following system of equations
\begin{equation}
\begin{cases}
a_1^2+a_2^2+a_3^2+a_4^2-2a_1a_3-2a_2a_4  =  8\\
2a_1a_2+2a_2a_3+2a_3a_4-6a_1a_4  =  -\frac{29}{2}\\
4a_1a_3+4a_2a_4  =  10\\
8a_1a_4=-\frac{5}{2}
\end{cases}
\end{equation}
Simplify (5.1) we have the following system
\begin{equation}
\begin{cases}
a_1^2+a_2^2+a_3^2+a_4^2  =  13\\
a_1a_2+a_2a_3+a_3a_4  = -\frac{131}{16}  \\
a_1a_3+a_2a_4 =  \frac{5}{2}\\
a_1a_4=-\frac{5}{16}
\end{cases}
\end{equation}
Solving for this system of equations yields 8 solutions. One of the numerical solutions is
\begin{equation}
\begin{cases}
a_1=2.6064\\
a_2=-2.3381\\
a_3=0.8516\\
a_4=-0.1199
\end{cases}
\end{equation}
We verify the condition for $\phi_4$ to be in $L_2(\mathbb{R})$
\begin{align*}
a_1^2+a_2^2+a_3^2+a_4^2=13 < 32=\frac{2^{2 \cdot 4 -1}}{4}
\end{align*}
Thus $\phi_4$ is indeed in $L_2(\mathbb{R})$. Now we will attempt to construct $\phi_4$ explicitly. It has the mask
\begin{align*}
P_4(z)S_4(z)&=\left(\frac{1+z}{2}\right)^4 (a_1z+a_2z^2+a_3z^3+a_4z^4)\\
&=0.1629 z + 0.5055 z^2 + 0.4461 z^3 - 0.0198 z^4 - 0.1323 z^5 + 0.0218 z^6 \\
&+ 0.0233 z^7 - 0.0075 z^8
\end{align*}
Hence, we have the refinement equation
\begin{equation}
\begin{split}
\phi_4(x)&=0.3258 \phi_4(2x-1) + 1.011 \phi_4(2x-2) + 0.8922 \phi_4(2x-3) - 0.0396 \phi_4(2x-4)\\
&- 0.2646 \phi_4(2x-5) + 0.0436 \phi_4(2x-6)+0.0466\phi_4(2x-7)-0.015\phi_4(2x-8)
\end{split}
\end{equation}

The above examples demonstrate the efficiency of the computation of the Daubechies scaling functions by using expression \eqref{0-4}.

\end{document}